# Measurement to Error Stability: a Notion of Partial Detectability for Nonlinear Systems

Brian P. Ingalls[*]
Control and Dynamical Systems, California Institute of Technology, CA
`ingalls@cds.caltech.edu`

Eduardo D. Sontag[†]
Dept. of Mathematics, Rutgers University, NJ
`sontag@math.rutgers.edu`

Yuan Wang[‡]
Dept. of Mathematics, Florida Atlantic University, FL
`ywang@math.fau.edu`

**Abstract**

In previous work the notion of input to state stability (ISS) has been generalized to systems with outputs, yielding a number of useful concepts. When considering a system whose output is to be kept small (i.e. an error output), the notion of input to output stability (IOS) arises. Alternatively, when considering a system whose output is meant to provide information about the state (i.e. a measurement output), one arrives at the detectability notion of output to state stability (OSS). Combining these concepts, one may consider a system with *two* outputs, an error and a measurement. This leads naturally to a notion of *partial detectability* we call measurement to error stability (MES). This property characterizes systems in which the error signal is detectable through the measurement signal.

This paper provides a partial Lyapunov characterization of the MES property. A closely related property of *stability in three measures* (SIT) is introduced, which characterizes systems for which the error decays whenever it dominates the measurement. The SIT property is shown to imply MES, and the two are shown to be equivalent under an additional boundedness assumption. A nonsmooth Lyapunov characterization of the SIT property is provided, which yields the partial characterization of MES. The analysis is carried out on systems described by differential inclusions – implicitly incorporating a disturbance input with compact value-set.

## 1 Introduction

The notion of *input to state stability* (ISS), introduced in [24], provides a theoretical framework in which to formulate questions of robustness with respect to inputs (seen as disturbances) acting on a system. An ISS system is, roughly, one which has a "finite nonlinear gain" with respect to inputs and whose transient behavior can be bounded in terms of the size of the initial state and inputs; the precise definition is in terms of $\mathcal{K}$-function gains. The theory of ISS systems now forms an integral part of several texts ([4, 10, 12, 15, 16, 23]) as well as expository and research articles (see e.g. [11, 13, 18, 20, 22, 26, 31, 33]).

In light of the duality between input/state and state/output behaviour which is common in control theory, it is natural to ask whether an ISS-like notion of output to state stability can be formulated. This concept, called OSS, is the subject of [14, 27, 28]. The definition given is precisely the same as that of ISS with outputs in the place of inputs. In the case of linear systems this property is equivalent to *detectability*. (When applied to nonlinear systems, OSS is more properly described as zero-detectability).

[*]research completed while the author was at Rutgers University, supported in part by US Air Force Grant F49620-98-1-0242
[†]Supported in part by US Air Force Grant F49620-01-1-0063
[‡]Supported in part by NSF Grant DMS-0072620

The paper [28] contains a discussion of various definitions of detectability for nonlinear systems which have appeared in the literature. Several of these definitions are given in terms of the existence of a Lyapunov or "storage" function for the system. The main result of [27] is the fact that the OSS property is equivalent to the existence of an appropriate Lyapunov function. These papers also contain a discussion of a generalized notion in which both inputs and outputs are considered (input-output to state stability, or IOSS). This property is addressed more completely in [14] where a Lyapunov characterization is provided and the construction of nonlinear observers is discussed.

This work addresses a generalization of the OSS property to a notion of "partial detectability". When discussing systems with outputs, the output signal typically plays one of two roles. A common situation is when the outputs are considered as *measurements*. Here, one supposes that knowledge of the whole state is not available, but rather that only partial knowledge of the state can be used. (Most commonly the output map is a projection, which corresponds simply to the ability to measure some, but not all, of the components of the state. More generally, one may only have access to some function of the state variables – e.g. the sum of two components – and so we allow for more general output mappings in the theory). This is the role of the output in OSS, and in the theory of detectability and observers in general.

A second role for outputs occurs when the goal of the control design is not to regulate the behaviour of the entire state, but rather only to regulate the output signal. The theory of *output regulation* addresses precisely this situation (see e.g. [9]). In the case of systems with no inputs, the problem of *stability* of a subset of the state variables (i.e. stability of an output signal which is a projection) has been addressed in the ordinary differential equations literature under the name "partial stability" [34]. Within the ISS framework, the notion of stability of the output signal has been described by *input to output stability* (IOS) [6, 29, 30].

Consider now the case in which both the above situations occur. That is, there are *two* output signals, one which is measured, and the other which must be regulated. A special case of this situation has been addressed in the output regulation theory, under the name "error feedback". This theory formulates the question of regulating an output of the system (the error) with knowledge of that output only. The more general case is when there are two distinct channels playing these two roles. In this paper we generalize the notion of OSS to this situation by introducing the concept of *measurement to error stability* (MES), which can be viewed as a notion of *partial detectability* through the measurement channel.

In this paper we will present a partial Lyapunov characterization of the MES property. This will be accomplished by first comparing the MES property to a notion of output stability relative to a set. This notion, which will be called *stability in three measures* (SIT) (cf. [17]) will be characterized by the existence of a lower semicontinuous Lyapunov function. It will be shown that the SIT property implies MES, and that the converse holds under an additional boundedness assumption.

All stability notions discussed in this paper are defined "robustly" with respect to disturbances. Disturbances are incorporated implicitly into the model by describing the dynamics of the system by a differential inclusion.

## 2 Basic Definitions and Notations

We consider the differential inclusion

$$\dot{x}(t) \in F(x(t)) \tag{1}$$

with two output maps

$$y(t) = h(x(t)), \quad w(t) = g(x(t)),$$

and a map $\omega : \mathbb{R}^n \to \mathbb{R}_{\geq 0}$. We take the state $x \in \mathbb{R}^n$. We assume that the set-valued map $F$ from $\mathbb{R}^n$ to subsets of $\mathbb{R}^n$ is locally Lipschitz (precise definitions to follow) with nonempty compact values. In addition, we assume that the differential inclusion (1) is forward complete. We assume that the output maps $h : \mathbb{R}^n \to \mathbb{R}^{p_y}$ and $g : \mathbb{R}^n \to \mathbb{R}^{p_w}$ are locally Lipschitz. The map $\omega$ is assumed to be continuous and proper; it will be used as a measurement of the magnitude of the state vector. We will denote $|\cdot|_\omega := \omega(\cdot)$. The use of $|\cdot|_\omega$ allows a framework which includes the Euclidean norm, distance to a compact set, and more general measures of the magnitude of the state.

**Remark 2.1** We use the setting of a differential inclusion as a generalization of the perturbed differential equation

$$\dot{x} = f(x, d) \qquad (2)$$

where $f$ is locally Lipschitz and the inputs $d(\cdot)$, thought of as disturbances, take values in some compact set $D$. The setup provided by (1) includes this case by choosing $F(x) := \{f(x, d) : d \in D\}$. □

The Euclidean norm in a space $\mathbb{R}^k$ is denoted simply by $|\cdot|$. If $z$ is a function defined on a real interval containing $[0, t]$, $\|z\|_{[0,t]}$ is the sup norm of the restriction of $z$ to $[0, t]$, that is $\|z\|_{[0,t]} =$ ess sup $\{|z(t)| : t \in [0, t]\}$. For each $p \in \mathbb{R}^n$ and $r \geq 0$ let $B(p, r) := \{x \in \mathbb{R}^n : |x - p| \leq r\}$, the ball of radius $r$ centered at $p$. Let $\mathcal{B}$ denote the unit ball $B(0, 1)$.

To formulate the statement that a nonsmooth function decreases in an appropriate manner, we will make use of the notion of the *viscosity subgradient* (cf. [1]).

**Definition 2.2** A vector $\zeta \in \mathbb{R}^n$ is a *viscosity subgradient* of the function $V : \mathbb{R}^n \to \mathbb{R}$ at $\xi \in \mathbb{R}^n$ if there exists a function $g : \mathbb{R}^n \to \mathbb{R}$ satisfying $\lim_{h \to 0} \frac{g(h)}{|h|} = 0$ and a neighbourhood $\mathcal{O} \subset \mathbb{R}^n$ of the origin so that

$$V(\xi + h) - V(\xi) - \zeta \cdot h \geq g(h)$$

for all $h \in \mathcal{O}$.

The (possibly empty) set of viscosity subgradients of $V$ at $\xi$ is called the *viscosity subdifferential* and is denoted $\partial_D V(\xi)$. We remark that if $V$ is differentiable at $\xi$, then $\partial_D V(\xi) = \{\nabla V(\xi)\}$.

A function $\gamma : \mathbb{R}_{\geq 0} \to \mathbb{R}_{\geq 0}$ is *of class* $\mathcal{K}$ (denoted $\gamma \in \mathcal{K}$) if it is continuous, positive definite, and strictly increasing; and is of class $\mathcal{K}_\infty$ if in addition it is unbounded. A function $\beta : \mathbb{R}_{\geq 0} \times \mathbb{R}_{\geq 0} \to \mathbb{R}_{\geq 0}$ is *of class* $\mathcal{KL}$ if for each fixed $t \geq 0$, $\beta(\cdot, t)$ is of class $\mathcal{K}$ and for each fixed $s \geq 0$, $\beta(s, t)$ decreases to zero as $t \to \infty$.

We next cite two results on nonlinear gain functions. The first is a small-gain lemma which is a special case of the main result in [7].

**Lemma 2.3** Suppose given a $\mathcal{KL}$ function $\beta$, a function $R : \mathbb{R}_{\geq 0} \times \mathbb{R}_{\geq 0} \to \mathbb{R}_{\geq 0}$, a $\mathcal{K}$ function $\gamma$ for which $\gamma(r) < r$ if $r > 0$. Then if a system as in (1) and a time $t_1 \geq 0$ satisfy the following:

i: for each $0 \leq t_0 \leq t \leq t_1$

$$|y(t)| \leq \max\{\beta(|x(t_0)|_\omega, t - t_0), \gamma(\|y\|_{[t_0,t]})\},$$

ii: for each $0 \leq t_0 \leq t < t_1$, the state satisfies the reachability condition

$$|x(t)|_\omega \leq R(r, t - t_0),$$

if $|x(t_0)|_\omega \leq r$;

then there exists a $\mathcal{KL}$ function $\widetilde{\beta}$ which satisfies

$$|y(t)| \leq \widetilde{\beta}(|x(t_0)|_\omega, t - t_0)$$

for each trajectory of the system and for all $0 \leq t_0 \leq t < t_1$. □

The next proposition follows directly from the proof of Lemma 3.1 in [19].

**Proposition 2.4** For any given $\mathcal{KL}$ function $\beta$, there exist a family of mappings $\{T_r\}_{r \geq 0}$ with:

- for each fixed $r > 0$, $T_r : \mathbb{R}_{>0} \overset{\text{onto}}{\to} \mathbb{R}_{>0}$ is continuous and is strictly decreasing;
- for each fixed $\varepsilon > 0$, $T_r(\varepsilon)$ is strictly increasing as $r$ increases and $\lim_{r \to \infty} T_r(\varepsilon) = \infty$;

such that

$$\beta(s, t) \leq \varepsilon$$

for all $s \leq r$, all $t \geq T_r(\varepsilon)$. □

## 2.1 Differential Inclusions

We review some standard concepts from set-valued analysis (See e.g. [1, 2, 3]). The following statements apply to a map $F$ from $\mathbb{R}^n$ to subsets of $\mathbb{R}^n$.

**Definition 2.5** Let $0 < T \leq \infty$. A function $x : [0, T) \to \mathbb{R}^n$ is said to be a *solution of the differential inclusion* (1) if it is absolutely continuous and satisfies

$$\dot{x}(t) \in F(x(t)),$$

for almost every $t \in [0, T)$. A function $x : [0, T) \to \mathbb{R}^n$ is called a *maximal solution of the differential inclusion* (1) if it does not have an extension which is a solution. That is, either $T = \infty$ or there does not exist a solution $\widehat{x} : [0, T_+) \to \mathbb{R}^n$ with $T_+ > T$ so that $\widehat{x}(t) = x(t)$ for all $t \in [0, T)$.

**Definition 2.6** The differential inclusion (1) is said to be *forward complete* on $\mathbb{R}^n$ if every maximal solution is defined for all $t \geq 0$.

For each $C \subseteq \mathbb{R}^n$ we let $\mathbf{S}(C)$ denote the set of maximal solutions of (1) satisfying $x(0) \in C$ equipped with the topology of uniform convergence on compact intervals. If $C$ is a singleton $\{\xi\}$ we will use the shorthand $\mathbf{S}(\xi)$. We set $\mathbf{S} := \mathbf{S}(\mathbb{R}^n)$, the set of all maximal solutions. Given a trajectory $x(\cdot) \in \mathbf{S}(\xi)$ for some $\xi \in \mathbb{R}^n$, we denote

$$y(t) = h(x(t)) \qquad w(t) = g(x(t)),$$

for all $t \geq 0$.

**Definition 2.7** Let $\mathcal{O}$ be an open subset of $\mathbb{R}^n$. The set-valued map $F$ is said to be *locally Lipschitz* on $\mathcal{O}$ if, for each $\xi \in \mathcal{O}$, there exists a neighbourhood $U \subset \mathcal{O}$ of $\xi$ and an $L > 0$ so that for any $\eta, \zeta$ in $U$,

$$F(\eta) \subseteq F(\zeta) + L |\eta - \zeta| \mathcal{B}.$$

An immediate consequence of the definition of a locally Lipschitz set-valued map is the following.

**Lemma 2.8** Suppose the set-valued map $F$ is locally Lipschitz on an open subset $\mathcal{O}$ of $\mathbb{R}^n$. Then, for any compact set $K \subset \mathcal{O}$, there exists some $L_K > 0$ so that for any $\eta, \zeta$ in $K$,

$$F(\eta) \subseteq F(\zeta) + L_K |\eta - \zeta| \mathcal{B}.$$

$\square$

We will use the notation $\mathcal{R}_{\leq T}(C)$ for the reachable set in time $T$ starting in the set $C$ for the differential inclusion (1). That is, for each $T > 0$ and $C \subseteq \mathbb{R}^n$,

$$\mathcal{R}_{\leq T}(C) := \{\eta \in \mathbb{R}^n \ : \ \eta = x(t) \ \text{ for some } \ x(\cdot) \in \mathbf{S}(C) \ \text{ and } \ t \in [0, T]\}.$$

The next result follows immediately from Corollary 3.4 of [8] and Theorem 3, §7, of [3].

**Lemma 2.9** Suppose the set-valued map $F$ from $\mathbb{R}^n$ to subsets of $\mathbb{R}^n$ is locally Lipschitz with nonempty compact values, and the differential inclusion as in (1) is forward complete. Then, for each $T \geq 0$ and each compact set $C \subset \mathbb{R}^n$ the set $\mathcal{R}_{\leq T}(C)$ is bounded. $\square$

The following generalization of Gronwall's Lemma will be needed. This is a special case of Lemma 8.3 in [2].

**Lemma 2.10** Suppose the set-valued map $G$ defined on $\mathbb{R}^n$ has closed nonempty values and is *globally* Lipschitz with constant $L$. Let $T > 0$ be given. Then for any solution $x(\cdot)$ of

$$\dot{x}(t) \in G(x(t)) \tag{3}$$

defined for $t \in [0, T]$ and any $p \in \mathbb{R}^n$, there is a solution $z_p(\cdot)$ of (3) defined on $[0, T]$ which has $z_p(0) = p$ and satisfies

$$|x(t) - z_p(t)| \leq |x(0) - p| \, e^{Lt} \qquad \forall t \in [0, T].$$

$\square$

We next make the straightforward observation that the result of the previous Lemma is valid for locally Lipschitz set-valued maps provided we restrict to compact sets.

**Lemma 2.11** Suppose given a system as in (1). Let a compact $C \subset \mathbb{R}^n$ and $T > 0$ be given. Then there exists $L > 0$ such that for any solution $x(\cdot)$ of (1) defined on $[0, T]$ which satisfies $x(0) \in C$ and any $p \in C$, there is a solution $z_p(\cdot)$ of (1) defined on $[0, T]$ which has $z_p(0) = p$ and satisfies

$$|x(t) - z_p(t)| \leq |x(0) - p| \, e^{Lt} \qquad \forall t \in [0, T].$$

*Proof.* Let $\Phi : \mathbb{R}^n \to [0, 1]$ be a smooth function so that $\Phi(x) = 1$ for all $x \in \mathcal{R}_{\leq T}(C)$ and $\Phi(x) = 0$ for all $x \notin B(\mathcal{R}_{\leq T}(C), 1)$. Since $\mathcal{R}_{\leq T}(C)$ is bounded (Lemma 2.9), the set-valued function $\widetilde{F}$ defined by $\widetilde{F}(x) = \Phi(x) F(x)$ is globally Lipschitz, say with constant $L$. The result follows from Lemma 2.10 and the fact that $\widetilde{F}$ and $F$ agree on the set $\mathcal{R}_{\leq T}(C)$, which contains the trajectories of interest. ∎

## 3 Stability and Detectability Properties

The following definitions are given for a forward complete system with two output channels as in (1). The outputs $y$ and $w$ are considered as error and measurement signals, respectively.

Our primary motivation is the following notion.

**Definition 3.1** We say that the system (1) is *measurement to error stable* (MES) if there exist $\beta \in \mathcal{KL}$ and $\gamma \in \mathcal{K}$ so that

$$|y(t)| \leq \max\{\beta(|x(0)|_\omega, t), \gamma(\|w\|_{[0,t]})\}$$

for each $x(\cdot) \in \mathbf{S}$, and all $t \geq 0$.

In the investigation of the MES property, the following notion of relative stability of the error will be useful. This is a notion of output stability which is applicable to systems with a single output $y$.

**Definition 3.2** Given a closed subset $D$ of the state space $\mathbb{R}^n$, we say that the system (1) is *relatively error stable (RES) with respect to $D$* if there exists $\beta \in \mathcal{KL}$ so that for any solution $x(\cdot) \in \mathbf{S}$, if there exists $t_1 > 0$ so that $x(t) \notin D$ for all $t \in [0, t_1]$, then

$$|y(t)| \leq \beta(|x(0)|_\omega, t) \qquad \forall t \in [0, t_1].$$

A special case of this property occurs for a system with two outputs when the set $D$ is defined by an inequality involving the two output maps, as follows.

**Definition 3.3** Let $\rho \in \mathcal{K}$. We say that the system (1) satisfies the *stability in three measures* (SIT) property (with gain $\rho$) if there exists $\beta \in \mathcal{KL}$ so that for any solution $x(\cdot) \in \mathbf{S}$, if there exists $t_1 > 0$ so that $|y(t)| > \rho(|w(t)|)$ for all $t \in [0, t_1]$, then

$$|y(t)| \leq \beta(|x(0)|_\omega, t) \qquad \forall t \in [0, t_1].$$

It is immediate that SIT is equivalent to relative error stability with respect to the set $D := \{\xi \in \mathbb{R}^n : |h(\xi)| \leq \rho(|g(\xi)|)\}$.

The following relative stability properties will also be considered.

**Definition 3.4** We say the system (1) satisfies the *relative measurement to error bounded property* (RMEB) if there exist $\mathcal{K}$ functions $\rho_1$, $\sigma_1$, and $\sigma_2$ so that for any solution $x(\cdot) \in \mathbf{S}$, if there exists $t_1 > 0$ so that $|y(t)| > \rho_1(|w(t)|)$ for all $t \in [0, t_1]$, then

$$|y(t)| \leq \max\{\sigma_1(|h(x(0))|), \sigma_2(\|w\|_{[0,t]})\} \qquad \forall t \in [0, t_1]. \tag{4}$$

**Definition 3.5** We say the system (1) satisfies the *relative error bounded property* (REB) if there exist $\mathcal{K}$ functions $\rho_2$ and $\sigma$ so that for any solution $x(\cdot) \in \mathbf{S}$, if there exists $t_1 > 0$ so that $|y(t)| > \rho_2(|w(t)|)$ for all $t \in [0, t_1]$, then

$$|y(t)| \leq \sigma(|h(x(0))|) \qquad \forall t \in [0, t_1].$$

We begin with the straightforward observation that the REB and RMEB properties are equivalent.

**Lemma 3.6** The system (1) satisfies the RMEB property if and only if it satisfies the REB property.

*Proof.* One implication is immediate.

Suppose the system (1) satisfies the RMEB property with $\rho_1$, $\sigma_1$, $\sigma_2 \in \mathcal{K}$. Define $\rho_2 \in \mathcal{K}$ by $\rho_2(r) = \max\{\rho_1(r), \sigma_2(r)\}$ for all $r \geq 0$. We will show that the system satisfies the REB property for this $\rho_2$ and with $\sigma = \sigma_1$.

Suppose given a trajectory $x(\cdot) \in \mathbf{S}$ which satisfies $|y(t)| > \rho_2(|w(t)|)$ on some interval $[0, t_1]$. Since $\rho_2 \geq \rho_1$, the RMEB bound (4) gives, for each $t \in [0, t_1]$,

$$\begin{aligned} |y(t)| &\leq \max\{\sigma_1(|h(x(0))|), \sigma_2(\|w\|_{[0,t]})\} \\ &\leq \max\{\sigma_1(|h(x(0))|), \sigma_2(\|w\|_{[0,t_1]})\}. \end{aligned}$$

Taking the supremum over $t \in [0, t_1]$, we have

$$\|y\|_{[0,t_1]} \leq \max\{\sigma_1(|h(x(0))|), \sigma_2(\|w\|_{[0,t_1]})\}. \tag{5}$$

Since $\rho_2 \geq \sigma_2$, we have

$$|y(t)| > \sigma_2(|w(t)|)$$

for all $t \in [0, t_1]$. Taking suprema over such $t$, we find

$$\|y\|_{[0,t_1]} > \sigma_2(\|w\|_{[0,t_1]}).$$

Finally, (5) gives

$$\|y(t)\|_{[0,t_1]} \leq \sigma_1(|h(x(0))|).$$

Thus the system satisfies the REB property. ∎

In the next section we provide a Lyapunov characterization for the relative error stability property. Since error stability relative to a set $D$ involves a condition only on $\mathbb{R}^n \setminus D$, one would expect a necessary and sufficient Lyapunov condition to be the existence of a function which decays appropriately along trajectories in $\mathbb{R}^n \setminus D$. This is the case. Unfortunately, the "natural" method of construction leads to a function which is only lower semicontinuous.

## 4 Lyapunov Functions

We give definitions of the appropriate Lyapunov functions.

**Definition 4.1** Given an open set $E \subseteq \mathbb{R}^n$, we say that a lower semicontinuous function $V : \mathbb{R}^n \to \mathbb{R}_{\geq 0}$ is a *lower semicontinuous RES-Lyapunov function* for system (1) on $E$ if

- there exist $\alpha_1$, $\alpha_2 \in \mathcal{K}_\infty$ so that

$$\alpha_1(|h(\xi)|) \leq V(\xi) \leq \alpha_2(|\xi|_\omega), \qquad \forall \xi \in E, \tag{6}$$

- there exists $\alpha_3 : \mathbb{R}_{\geq 0} \to \mathbb{R}_{\geq 0}$ continuous positive definite so that for each $\xi \in E$,

$$\zeta \cdot v \leq -\alpha_3(V(\xi)) \qquad \forall \zeta \in \partial_D V(\xi),\ \forall v \in F(\xi). \tag{7}$$

We say that $V$ is a *lower semicontinuous exponential decay RES-Lyapunov function* for system (1) on $E$ if in addition (7) holds with $\alpha_3(r) = r$.

We specialize the above definitions for the notion of stability in three measures as follows.

**Definition 4.2** Let $\rho \in \mathcal{K}$. We say that a lower semicontinuous function $V : \mathbb{R}^n \to \mathbb{R}_{\geq 0}$ is a *lower semicontinuous SIT-Lyapunov function* for system (1) with gain $\rho$ if

- there exist $\alpha_1$, $\alpha_2 \in \mathcal{K}_\infty$ so that

$$\alpha_1(|h(\xi)|) \leq V(\xi) \leq \alpha_2(|\xi|_\omega), \qquad \forall \xi \text{ so that } |h(\xi)| > \rho(|g(\xi)|),$$

- there exists $\alpha_3 : \mathbb{R}_{\geq 0} \to \mathbb{R}_{\geq 0}$ continuous positive definite so that for each $\xi$ so that $|h(\xi)| > \rho(|g(\xi)|)$,

$$\zeta \cdot v \leq -\alpha_3(V(\xi)) \qquad \forall \zeta \in \partial_D V(\xi),\ \forall v \in F(\xi). \tag{8}$$

We say that $V$ is a *lower semicontinuous exponential decay SIT-Lyapunov function* for system (1) with gain $\rho$ if in addition (8) holds with $\alpha_3(r) = r$.

We next remark that the decrease statements (7) and (8) can be written equivalently in an integral formulation, using the following standard result (a minor extension of Theorem 4.6.3 in [1], see e.g. [21] for details).

**Proposition 4.3** Suppose given a forward complete system

$$\dot{x} \in F(x)$$

where $F$ is locally Lipschitz and takes nonempty compact values. Let a lower semicontinuous function $V : \mathbb{R}^n \to \mathbb{R}_{\geq 0}$ and a locally Lipschitz function $w : \mathbb{R}^n \to \mathbb{R}$ be given. The following are equivalent:

1. For each $\xi \in \mathbb{R}^n$,

$$\zeta \cdot v \leq w(\xi) \qquad \forall \zeta \in \partial_D V(\xi), \ \forall v \in F(\xi).$$

2. For each $\xi \in \mathbb{R}^n$, each solution $x(\cdot) \in \mathbf{S}(\xi)$ verifies

$$V(x(t)) - V(\xi) \leq \int_0^t w(x(s))\, ds$$

for any $t \geq 0$.

□

Making use of this result, the decrease statements (7) and (8) above can be written equivalently as (after possibly replacing $\alpha_3$ by a locally Lipschitz function dominated by the original $\alpha_3$)

$$V(x(t)) - V(x(0)) \leq -\int_0^t \alpha_3(V(x(s)))\, ds, \tag{9}$$

for all $x(\cdot) \in \mathbf{S}$ which remain in the appropriate set on the interval $[0, t]$. This alternative formulation will be used below.

The Lyapunov characterizations are as follows.

**Theorem 1** *Let a system of the form (1) and a closed set $D \subset \mathbb{R}^n$ be given. Let $E = \mathbb{R}^n \setminus D$. The following are equivalent.*

1. *The system is relatively error stable with respect to $D$.*

2. *The system admits a lower semicontinuous RES-Lyapunov function on $E$.*

3. *The system admits a lower semicontinuous exponential decay RES-Lyapunov function on $E$.*

The implication (3) ⇒ (2) is immediate. The others will be shown in Section 6.

**Corollary 4.4** *Let a system of the form (1) and a function $\rho \in \mathcal{K}$ be given. The following are equivalent.*

1. *The system satisfies the SIT property with gain $\rho$.*

2. *The system admits a lower semicontinuous SIT-Lyapunov function with gain $\rho$.*

3. *The system admits a lower semicontinuous exponential decay SIT-Lyapunov function with gain $\rho$.*

□

The corollary follows immediately by setting $D = \{\xi \in \mathbb{R}^n\ :\ |h(\xi)| \leq \rho(|g(\xi)|)\}$.

# 5 Relationships between Notions

Having given a characterization of the SIT property, we now indicate how this notion is related to measurement to error stability. The following will be shown.

**Lemma 5.1** If the system (1) satisfies the MES property, then it satisfies the SIT property.

**Lemma 5.2** If the system (1) satisfies the RMEB property and the SIT property, then it satisfies the MES property.

In addition, we provide an example to show that the converse of Lemma 5.1 does not hold in general.

The following partial characterization of MES is an immediate consequence of Corollary 4.4 and the two preceding lemmas.

**Corollary 5.3** *If the system (1) satisfies MES, then it admits a lower semicontinuous exponential decay SIT-Lyapunov function. If the system satisfies the RMEB property and admits a lower semicontinuous SIT-Lyapunov function, then it satisfies MES.*  □

## 5.1 Proofs

We first show that the MES property implies the SIT property.
*Proof.* (Lemma 5.1)

Assume that the system (1) satisfies the MES property with gains $\beta$ and $\gamma$. Let $\rho$ be any $\mathcal{K}_\infty$ function so that $\gamma(\rho^{-1}(s)) < s$ for all $s > 0$.

For each $x(\cdot) \in \mathbf{S}$, it follows that if $t_1 > 0$ is such that
$$|y(t)| > \rho(|w(t)|) \qquad \forall t \in [0, t_1],$$
then from the definition of MES,
$$\begin{aligned} |y(t)| &\leq& \max\{\beta(|x(0)|_\omega, t), \gamma(\|w\|_{[0,t]})\} \\ &\leq& \max\{\beta(|x(0)|_\omega, t), \gamma(\rho^{-1}(\|y\|_{[0,t]}))\} \end{aligned}$$
for any $t \in [0, t_1]$. Lemma 5.1 follows from an application of Lemma 2.3 with this $\beta$ and $t_1$, and with $\gamma(\rho^{-1}(r))$ in the place of $\gamma$. The existence of the map $R(\cdot, \cdot)$ follows from Lemma 2.9.  ∎

We next show that under the RMEB condition, SIT implies MES. We will show the following, which gives Lemma 5.2 immediately.

**Lemma 5.4** If the system (1) satisfies the RMEB property with gains $\rho_1$, $\sigma_1$, and $\sigma_2$ and also satisfies the SIT property with gains $\widetilde{\rho} \in \mathcal{K}$, $\beta \in \mathcal{KL}$, then there exists a time $0 \leq t_1 \leq \infty$ so that

- $|y(t)| \leq \beta(|x(0)|_\omega, t)$ for all $t \in [0, t_1)$, and
- $|y(t)| \leq \gamma(\|w\|_{[t_1,t]})$ for all $t \geq t_1$,

where $\gamma(r) = \max\{\rho_1(r), \widetilde{\rho}(r), \sigma_1(\rho_1(r)), \sigma_1(\widetilde{\rho}(r)), \sigma_2(r)\}$.

*Proof.* Suppose a system satisfies RMEB and SIT with gains as above. Define the $\mathcal{K}$ functions $\rho(r) := \max\{\widetilde{\rho}(r), \rho_1(r)\}$ and $\gamma(r) := \max\{\rho(r), \sigma_1(\rho(r)), \sigma_2(r)\}$.

Let $x(\cdot) \in \mathbf{S}$ be given. Define
$$t_1 := \inf\{t \geq 0 \,:\, |y(t)| \leq \rho(|w(t)|)\},$$
with $t_1 = \infty$ if the inequality never holds. Let $t \geq 0$. We consider three possibilities (see Figure 1).
i) If $t \in [0, t_1)$, it follows that $|y(s)| > \rho(|w(s)|) \geq \widetilde{\rho}(|w(s)|)$ for all $s \in [0, t]$, so the SIT property gives
$$|y(t)| \leq \beta(|x(0)|_\omega, t).$$
ii) If $t \geq t_1$ and $|y(t)| \leq \rho(|w(t)|)$, then it is immediate that
$$|y(t)| \leq \rho(\|w\|_{[t_1,t]}) \leq \gamma(\|w\|_{[t_1,t]}).$$

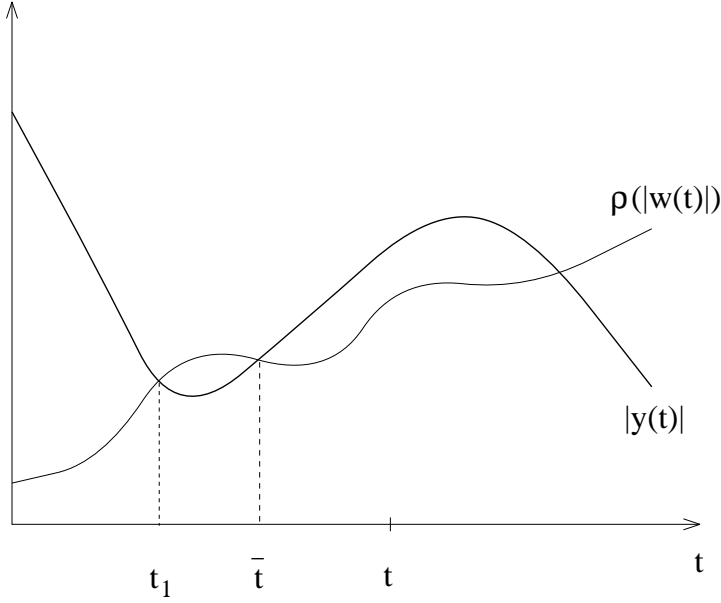

Figure 1: Proof of Lemma 5.4

iii) The last possibility is $t \geq t_1$ and $|y(t)| > \rho(|w(t)|)$. If this is the case, let
$$\bar{t} := \sup\{s \leq t \ : \ |y(s)| \leq \rho(|w(s)|)\}.$$
Note that there exist such $s$, since $y(t_1) \leq \rho(w(t_1))$ and $t \geq t_1$ (so $\bar{t} \in [t_1, t]$). Also, $|y(\bar{t})| = \rho(|w(\bar{t})|)$. Then, for each $0 < \varepsilon < t - \bar{t}$, we have
$$|y(s)| > \rho(|w(s)|) \geq \rho_1(|w(s)|) \qquad \forall s \in [\bar{t}+\varepsilon, t],$$
so the RMEB property gives
$$|y(t)| \ \leq \ \max\{\sigma_1(|y(\bar{t}+\varepsilon)|), \sigma_2(\|w\|_{[\bar{t}+\varepsilon,t]})\}.$$
As this holds for all $\varepsilon > 0$ sufficiently small, it follows by continuity that in the limit as $\varepsilon \to 0$,
$$\begin{aligned} |y(t)| &\leq& \max\{\sigma_1(|y(\bar{t})|), \sigma_2(\|w\|_{[\bar{t},t]})\} \\ &=& \max\{\sigma_1(\rho(|w(\bar{t})|)), \sigma_2(\|w\|_{[\bar{t},t]})\} \\ &\leq& \max\{\sigma_1(\rho(\|w\|_{[t_1,t]})), \sigma_2(\|w\|_{[t_1,t]})\} \\ &\leq& \gamma(\|w\|_{[t_1,t]}). \end{aligned}$$

■

To complement these results, we next exhibit an example showing the SIT property alone does not imply MES.

**Example 5.5** Consider the system evolving in $\mathbb{R}^2$ defined by the following differential equations
$$\begin{aligned} \dot{x}_1 &=& x_1 + x_2, \\ \dot{x}_2 &=& -x_1 + x_2, \end{aligned}$$

with outputs
$$h(x_1, x_2) = x_1 \qquad g(x_1, x_2) = 1.$$
As usual, denote $y(t) := h(x_1(t), x_2(t))$ and $w(t) := g(x_1(t), x_2(t)) \equiv 1$. Take $|\cdot|_\omega = |\cdot|$.
Define $\rho(s) := s$ and $\beta(s,t) := s e^\pi (e^{\pi - t})$.

For initial condition $(x_1(0), x_2(0)) = (\xi_1, \xi_2)$, the solution to the system is:
$$\begin{aligned} x_1(t) &= e^t(\xi_1 \cos t + \xi_2 \sin t) \\ x_2(t) &= e^t(-\xi_1 \sin t + \xi_2 \cos t). \end{aligned}$$

We next verify the SIT property for this system. Let an initial condition $\xi = (\xi_1, \xi_2)$ and a time $t_1 \geq 0$ be so that the trajectory $x(\cdot) = (x_1(\cdot), x_2(\cdot))$ starting at $\xi$ has $|y(t)| > \rho(|w(t)|)$ for $t \in [0, t_1]$. Let $M := |\xi|$. Then $|y(t)| \leq |x(t)| = Me^t$ for each $t \geq 0$.

Set $t_2 \in [0, \pi)$ so that $y(t_2) = 0$, i.e.
$$t_2 := \begin{cases} \tan^{-1}(\frac{\xi_1}{\xi_2}) & \text{if } \frac{\xi_1}{\xi_2} \geq 0, \xi_2 \neq 0 \\ \pi + \tan^{-1}(\frac{\xi_1}{\xi_2}) & \text{if } \frac{\xi_1}{\xi_2} < 0, \xi_2 \neq 0 \\ \pi/2 & \text{if } \xi_2 = 0. \end{cases}$$

Since $y(t_2) = 0$, we have $|y(t_2)| < \rho(|w(t_2)|)$, which implies $t_1 \in [0, t_2)$.

We note that $|y(t)| \leq Me^{t_2}$ for all $t \in [0, t_1]$. Thus, for each $t \in [0, t_1]$, we have, as $t \leq t_1 < t_2 < \pi$,
$$\begin{aligned} |y(t)| &\leq Me^{t_2} \\ &\leq Me^{t_2-\pi}e^\pi(e^{\pi-t}) \\ &= \beta(Me^{t_2-\pi}, t) \\ &\leq \beta(M, t) \\ &= \beta(|x(0)|, t). \end{aligned}$$

Thus the SIT property holds.

However it is clear that the MES property does not hold, since the oscillations of the error $y(\cdot)$ grow without bound, and so cannot verify the MES bound for any gains $\beta$ and $\gamma$.

**Remark 5.6** This example indicates the existence of a large class of systems which satisfy the SIT property. Suppose a forward complete system satisfies the growth condition
$$|x(t)| \leq \sigma(|x(0)|)\sigma(t) \qquad \forall t \geq 0$$
for some $\sigma \in \mathcal{K}$. Suppose further that for some $\rho \in \mathcal{K}$, $|y(t)|$ is never greater than $\rho(|w(t)|)$ for longer than $T$ time units. Assuming $h$ is continuous, and $h(0) = 0$, there is an $\alpha \in \mathcal{K}$ so that $|h(\xi)| \leq \alpha(|\xi|)$ for all $\xi \in \mathbb{R}^n$. Then, on any interval $[0, t_1]$ on which $|y(t)| > \rho(|w(t)|)$, it follows that $t_1 \leq T$, and so
$$|y(t)| \leq \alpha(|x(t)|) \leq \alpha(\sigma(|x(0)|)\sigma(T)) \qquad \forall t \in [0, t_1].$$
Thus the SIT property holds, e.g. with $\beta(r,t) = \alpha(\sigma(r)\sigma(T))e^{T-t}$, irrespective of the dynamics. □

# 6 Proof of Theorem 1

We first present the proof of (2) ⇒ (1) in Theorem 1 (sufficiency).

## 6.1 Sufficiency

We will show that the existence of a lower semicontinuous RES-Lyapunov function implies relative error stability. We will make use of the following comparison result.

**Lemma 6.1** Suppose given a locally Lipschitz positive definite function $\alpha : \mathbb{R}_{\geq 0} \to \mathbb{R}_{\geq 0}$. Let $0 < \widetilde{t} \leq \infty$, and $v : [0, \widetilde{t}) \to \mathbb{R}_{\geq 0}$ be any lower semicontinuous function which satisfies
$$v(t_2) \leq v(t_1) - \int_{t_1}^{t_2} \alpha(v(s))\, ds \qquad \forall\, 0 \leq t_1 \leq t_2 < \widetilde{t}. \tag{10}$$

Define $w(\cdot)$ to be the maximal solution of the initial value problem
$$\dot{w}(t) = -\alpha(w(t)), \qquad w(0) = v(0).$$

Then $w(t)$ is defined for all $t \geq 0$ and
$$v(t) \leq w(t) \qquad \forall t \in [0, \widetilde{t}).$$

*Proof.* (We follow the proof of Theorem III.4.1 in [5]). Let $v(\cdot)$, $w(\cdot)$ be as above for given $\alpha(\cdot)$ and $\widetilde{t}$. It is immediate that $w(t) \in [0, w(0)]$ for all $t$ for which $w(t)$ is defined, hence $w(\cdot)$ is defined for all $t \geq 0$. For each integer $n \geq 1$, let $w_n(\cdot)$ be the maximal solution of

$$\dot{w}_n(t) = -\alpha(w_n(t)) + \frac{1}{n}, \qquad w_n(0) = v(0). \tag{11}$$

Then for each $n$, $w_n(t) \in [0, w_n(0) + \frac{t}{n}]$ for all $t$ for which $w_n(t)$ is defined, hence each $w_n(\cdot)$ is defined for all $t \geq 0$. We will show that

$$v(t) \leq w_n(t) \qquad \forall t \in [0, \widetilde{t}), \tag{12}$$

for all $n \geq 1$. Indeed, suppose not. Then there exists $n \geq 1$ and $\tau \in [0, \widetilde{t})$ so that

$$v(\tau) > w_n(\tau).$$

Let $t_0 := \sup\{0 \leq t \leq \tau : v(t) \leq w_n(t)\}$. Then, as $v(\cdot)$ is lower semicontinuous and $w_n(\cdot)$ is continuous, $v(t_0) \leq w_n(t_0)$ (since $\{0 \leq t \leq \tau : v(t) \leq w_n(t)\}$ is closed). Moreover, since $v(\cdot)$ is non-increasing and $w_n(\cdot)$ is continuous, it must be the case that $v(t_0) = w_n(t_0)$ (since $v(t_0) < w_n(t_0)$ would imply $v(t_0 + \delta) < w_n(t_0 + \delta)$ for $\delta > 0$ small). Finally, since $\{0 \leq t \leq \tau : v(t) > w_n(t)\}$ is open, we have $v(t_0 + \varepsilon) > w_n(t_0 + \varepsilon)$ for $\varepsilon > 0$ sufficiently small. From (10) and Taylor's Theorem, we have that for $\varepsilon \in (0, \tau - t_0)$,

$$\begin{aligned}
v(t_0 + \varepsilon) &\leq v(t_0) - \int_{t_0}^{t_0+\varepsilon} \alpha(v(s)) \, ds \\
&= v(t_0) - \varepsilon \alpha(v(t_0)) + o(\varepsilon)
\end{aligned}$$

and from (11),

$$\begin{aligned}
w_n(t_0 + \varepsilon) &= w_n(t_0) + \int_{t_0}^{t_0+\varepsilon} -\alpha(w_n(s)) + \frac{1}{n} \, ds \\
&= w_n(t_0) - \varepsilon \alpha(w_n(t_0)) + \frac{\varepsilon}{n} + o(\varepsilon),
\end{aligned}$$

where $o(\cdot)$ signifies a function satisfying $\lim_{t \to 0} \frac{o(t)}{t} = 0$. Since $w_n(t_0) = v(t_0)$, we conclude that $v(t_0 + \varepsilon) \leq w_n(t_0 + \varepsilon)$ for $\varepsilon > 0$ sufficiently small, which is a contradiction. We conclude that $v(t) \leq w_n(t)$ for all $t \in [0, \widetilde{t})$ and for all $n \geq 1$.

We note that $w_n(t) \to w(t)$ uniformly on each finite time interval. Thus for any $T \in [0, \widetilde{t})$, as (12) holds for all $n$,

$$v(t) \leq \lim_{n \to \infty} w_n(t) = w(t) \qquad \forall t \in [0, T].$$

As $T$ was arbitrary, we conclude that $v(t) \leq w(t)$ for all $t \in [0, \widetilde{t})$. ∎

The following is an immediate consequence of Lemma 6.1 and ([19], Lemma 4.4).

**Lemma 6.2** Suppose given a locally Lipschitz positive definite function $\alpha : \mathbb{R}_{\geq 0} \to \mathbb{R}_{\geq 0}$. Then there exists $\beta \in \mathcal{KL}$ with the following property: For any $0 < \widetilde{t} \leq \infty$ and for any lower semicontinuous function $v : [0, \widetilde{t}) \to \mathbb{R}_{\geq 0}$ which satisfies

$$v(t_2) \leq v(t_1) - \int_{t_1}^{t_2} \alpha(v(s)) \, ds \qquad \forall \, 0 \leq t_1 \leq t_2 < \widetilde{t},$$

it follows that

$$v(t) \leq \beta(v(0), t) \qquad \forall t \in [0, \widetilde{t}).$$

*Proof.* Lemma 6.1 tells us that $v(t) \leq w(t)$ for all $t \in [0, \widetilde{t})$, where $w$ is the maximal solution of the initial value problem

$$\dot{w}(t) = -\alpha(w(t)), \qquad w(0) = v(0).$$

Lemma 4.4 of [19] provides the existence of a $\beta \in \mathcal{KL}$ so that

$$w(t) \leq \beta(w(0), t) \qquad \forall t \geq 0.$$

Since $w(0) = v(0)$, we conclude that

$$v(t) \leq \beta(v(0), t) \qquad \forall t \in [0, \widetilde{t}).$$

∎

Finally, we give the proof of (2) ⇒ (1) for Theorem 1 (sufficiency).

*Proof.* Let a system of the form (1) and a closed set $D \subset \mathbb{R}^n$ be given. Let $E = \mathbb{R}^n \backslash D$. Suppose there exists a lower semicontinuous RES-Lyapunov function $V$ for the system which satisfies (6) and (7) (and hence (9)) on $E$ with gains $\alpha_1, \alpha_2 \in \mathcal{K}_\infty$ and $\alpha_3$ continuous positive definite. We will verify that the system satisfies the relative error stability property with respect to $D$.

Let $x(\cdot) \in \mathbf{S}$, and suppose $t_1 > 0$ is such that $x(t) \in E$ for all $t \in [0, t_1]$. Then we have, from (9),

$$V(x(t)) - V(x(0)) \leq - \int_0^t \alpha_3(V(x(s))) \, ds \qquad \forall t \in [0, t_1].$$

From Lemma 6.2 there exists $\widetilde{\beta} \in \mathcal{KL}$ (depending only on $\alpha_3$) so that

$$V(x(t)) \leq \widetilde{\beta}(V(x(0)), t) \qquad \forall t \in [0, t_1).$$

With (6), this gives

$$|h(x(t))| \leq \alpha_1^{-1}(V(x(t))) \leq \alpha_1^{-1}(\widetilde{\beta}(V(x(0)), t)) \leq \alpha_1^{-1}(\widetilde{\beta}(\alpha_2(|x(0)|_\omega), t))$$

for $t \in [0, t_1)$. By continuity, this holds on $[0, t_1]$. This is the relative error stability property with $\beta(r, t) = \alpha_1^{-1}(\widetilde{\beta}(\alpha_2(r), t))$. ∎

We next give the proof of (1) ⇒ (3) in Theorem 1 (necessity of a lower semicontinuous exponential decay RES-Lyapunov function).

## 6.2 Necessity

For a system (1) which satisfies the relative error stability property with respect to a set $D$, we define, for each $\xi \in E := \mathbb{R}^n \backslash D$ and each trajectory $x(\cdot) \in \mathbf{S}(\xi)$, the first hitting time of $x(\cdot)$ into $D$ as follows.

**Definition 6.3** For each $x(\cdot) \in \mathbf{S}(E)$, we set $\theta(x(\cdot)) = \inf\{t \geq 0 : x(t) \in D\}$, with $\theta(x(\cdot)) = \infty$ if $x(t) \notin D$ for all $t \geq 0$.

**Lemma 6.4** The map $x(\cdot) \mapsto \theta(x(\cdot))$ is lower semicontinuous on the set $\mathbf{S}(E)$.

*Proof.* Suppose $x(\cdot) \in \mathbf{S}(E)$ and the sequence $x_k(\cdot)$ is such that $x_k(\cdot) \in \mathbf{S}(E)$ for each $k$ and $x_k(\cdot) \to x(\cdot)$ as elements of $\mathbf{S}$. Denote $\theta_k := \theta(x_k(\cdot))$, and $\theta_0 = \liminf_{k \to \infty} \theta_k$. We will show $\theta(x(\cdot)) \leq \theta_0$. Without loss of generality, we may assume $\theta_0 < \infty$.

Passing to a subsequence if necessary, we assume $\theta_k \to \theta_0$. Thus, there exists $K > 0$ such that $\theta_k \leq \theta_0 + 1$ for all $k \geq K$. Since $x_k(\cdot)$ converges uniformly to $x(\cdot)$ on the finite interval $[0, \theta_0 + 1]$, it follows from continuity of $x(\cdot)$ and each $x_k(\cdot)$ that

$$x(\theta_0) = \lim_{k \to \infty} x_k(\theta_k).$$

Since $D$ is closed and $x_k(\theta_k) \in D$ for each $k$, we have $x(\theta_0) \in D$. Hence $\theta(x(\cdot)) \leq \theta_0$. ∎

Note that in general the function $\theta$ is not upper semicontinuous, as indicated in Figure 2.

We now present the proof of (1) ⇒ (3) for Theorem 1 (necessity).

*Proof.* Suppose the system (1) satisfies the relative error stability property with respect to a set $D$. Let $E := \mathbb{R}^n \backslash D$. We will construct an exponential decay lower semicontinuous RES-Lyapunov function for (1) on $E$.

Note that by definition of $\theta$, for each $\xi \in E$ and each $x(\cdot) \in \mathbf{S}(\xi)$ it follows that

$$|y(t)| \leq \beta(|\xi|_\omega, t) \qquad \forall t \in [0, \theta(x(\cdot))). \tag{13}$$

Now, consider a choice of $\widetilde{\alpha} \in \mathcal{K}_\infty$, $\widetilde{\beta} \in \mathcal{KL}$, and a smooth strictly increasing function $l : \mathbb{R}_{\geq 0} \to \mathbb{R}_{>0}$ such that

$$\widetilde{\alpha}(\beta(s, t)) l(t) \leq \widetilde{\beta}(s, t) \tag{14}$$

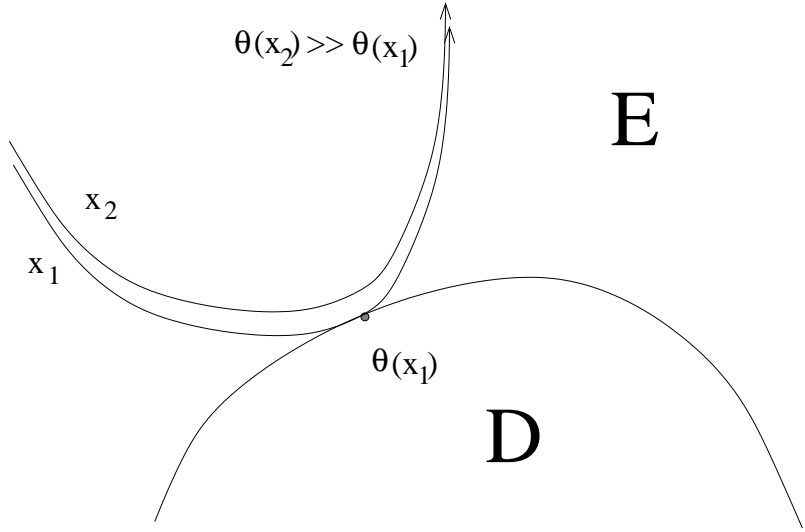

Figure 2: Proof of Theorem 1 – $\theta(\cdot)$ is not upper semicontinuous

for all $s$, $t \geq 0$. Such a set of functions always exists, e.g. $\widetilde{\alpha}(r) = r$, $\widetilde{\beta} = \beta$, and $l(t) \in (0,1)$ for all $t \geq 0$. Extend $l(\cdot)$ to $\mathbb{R}$ by setting $l(t) = l(0)$ for all $t < 0$.

For a given choice of $\widetilde{\alpha}$, $\widetilde{\beta}$ and $l(\cdot)$, we define a function $V : \mathbb{R}^n \to \mathbb{R}_{\geq 0}$ as follows. For each $\xi \in E$ we set

$$V(\xi) := \sup_{x(\cdot) \in \mathbf{S}(\xi)} \sup_{t \in [0, \theta(x(\cdot)))} \widetilde{\alpha}(|y(t)|)l(t),$$

and set $V(\xi) := 0$ for all $\xi \notin E$.

We note that

$$\widetilde{\alpha}(|h(\xi)|)l(0) \leq V(\xi) \leq \widetilde{\beta}(|\xi|_\omega, 0) \qquad \forall \xi \in E, \tag{15}$$

which follows immediately from

$$V(\xi) \geq \sup_{x(\cdot) \in \mathbf{S}(\xi)} \widetilde{\alpha}(|y(0)|)l(0) = \widetilde{\alpha}(|h(\xi)|)l(0) \qquad \forall \xi \in E,$$

and, from (13) and (14), for all $\xi \in E$,

$$\begin{aligned} V(\xi) &\leq \sup_{x(\cdot) \in \mathbf{S}(\xi)} \sup_{t \in [0, \theta(x(\cdot)))} \widetilde{\beta}(|\xi|_\omega, t) \\ &\leq \widetilde{\beta}(|\xi|_\omega, 0). \end{aligned} \tag{16}$$

**Lemma 6.5** The function $V$ is lower semicontinuous on $\mathbb{R}^n$.

*Proof.* It is immediate that $V$ is lower semicontinuous on $D$, since $V(\xi) = 0$ for all $\xi \in D$, and $V(\xi) \geq 0$ for all $\xi \in \mathbb{R}^n$.

Suppose $\xi \in E$, and $\xi_k \to \xi$ in $E$. We will show that

$$\liminf_{k \to \infty} V(\xi_k) \geq V(\xi).$$

Let $\varepsilon > 0$ be given. Pick a trajectory $\widehat{x}(\cdot) \in \mathbf{S}(\xi)$ for which

$$V(\xi) \leq \sup_{t \in [0, \theta(\widehat{x}(\cdot)))} \widetilde{\alpha}(|\widehat{y}(t)|)l(t) + \frac{\varepsilon}{3},$$

where $\widehat{y}(t) = h(\widehat{x}(t))$. Next choose $\tau < \theta(\widehat{x}(\cdot))$ so that

$$\sup_{t \in [0, \theta(\widehat{x}(\cdot)))} \widetilde{\alpha}(|\widehat{y}(t)|)l(t) \leq \max_{t \in [0, \tau]} \widetilde{\alpha}(|\widehat{y}(t)|)l(t) + \frac{\varepsilon}{3}.$$

Finally, choose any sequence $x_k(\cdot) \in \mathbf{S}$ so that each $x_k(\cdot) \in \mathbf{S}(\xi_k)$, and $x_k(\cdot) \to \widehat{x}(\cdot)$ uniformly on $[0,\tau]$ (such a sequence exists by Lemma 2.11). Denote $y_k(t) = h(x_k(t))$.

Since $x_k(\cdot) \to \widehat{x}(\cdot)$ uniformly on the interval $[0,\tau]$, there is some compact set $C \subset \mathbb{R}^n$ which contains the restriction to $[0,\tau]$ of all of the trajectories $\widehat{x}(\cdot)$ and $x_k(\cdot)$ for $k$ larger than some $K_1$. Then, since $\widetilde{\alpha}(|h(\cdot)|)$ is uniformly continuous on $C$, we can find $K_2 > K_1$ so that

$$|\widetilde{\alpha}(|h(\widehat{x}(t))|) - \widetilde{\alpha}(|h(x_k(t))|)| \leq \frac{\varepsilon}{3l(\tau)}$$

for each $k > K_2$ and each $t \in [0,\tau]$. Since $\tau < \theta(\widehat{x}(\cdot))$ and $x_k(\cdot) \to \widehat{x}(\cdot)$ uniformly on $[0,\tau]$, it follows that there exists $K_3 > K_2$ so that for $k > K_3$, $\theta(x_k(\cdot)) > \tau$. Finally, for each $k > K_3$, we have

$$\begin{aligned}
V(\xi_k) &= \sup_{x(\cdot) \in \mathbf{S}(\xi_k)} \sup_{t \in [0,\theta(x(\cdot)))} \widetilde{\alpha}(|y(t)|)l(t) \\
&\geq \sup_{t \in [0,\theta(x_k(\cdot)))} \widetilde{\alpha}(|y_k(t)|)l(t) \\
&\geq \max_{t \in [0,\tau]} \widetilde{\alpha}(|y_k(t)|)l(t) \\
&\geq \max_{t \in [0,\tau]} \widetilde{\alpha}(|\widehat{y}(t)|)l(t) - \frac{\varepsilon}{3} \\
&\geq \sup_{t \in [0,\theta(\widehat{x}(\cdot)))} \widetilde{\alpha}(|\widehat{y}(t)|)l(t) - \frac{2\varepsilon}{3} \\
&\geq V(\xi) - \varepsilon.
\end{aligned}$$

Hence

$$\liminf_{k \to \infty} V(\xi_k) \geq V(\xi) - \varepsilon.$$

Letting $\varepsilon$ tend to 0, we conclude that $V$ is lower semicontinuous at $\xi$. Since $\xi \in E$ was chosen arbitrarily, we have that $V$ is lower semicontinuous on $E$. ∎

We next consider the manner in which $V$ decreases along trajectories. Note that (15) and (16) give

$$V(\xi) = \sup_{x(\cdot) \in \mathbf{S}(\xi)} \sup_{t \in [0,\min\{\theta(x(\cdot)),T_\xi\})} \widetilde{\alpha}(|y(t)|)l(t).$$

where

$$T_\xi := \begin{cases} 0 & \text{if } |\xi|_\omega = 0 \\ T_{2|\xi|_\omega}\left(\widetilde{\alpha}\left(|h(\xi)/2|\right)l(0)\right) & \text{if } |\xi|_\omega > 0, |h(\xi)| > 0 \\ \infty & \text{if } |\xi|_\omega > 0, |h(\xi)| = 0, \end{cases}$$

for the function $T_r(\varepsilon)$ defined for $\widetilde{\beta}$ as in Proposition 2.4. (For consistency of notation, we interpret $\sup_{t \in [0,\min\{\theta(x(\cdot)),T_\xi\})} \widetilde{\alpha}(|y(t)|)l(t) = \widetilde{\alpha}(|y(0)|)l(0)$ if $T_\xi = 0$).

Let $\xi \in E$, and $x(\cdot) \in \mathbf{S}(\xi)$ so that $x(t) \in E$ for $t$ in some interval $[0,\overline{t}]$. For $\tau \in [0,\overline{t}]$ small, $x(\tau)$ satisfies $2|\xi|_\omega > |x(\tau)|_\omega$ and $\frac{|h(\xi)|}{2} < |h(x(\tau))|$, so the supremum in time in the expression for $V(x(\tau))$ may be taken over $[0,T_\xi]$. For such $\tau$, we find

$$\begin{aligned}
V(x(\tau)) &= \sup_{z(\cdot) \in \mathbf{S}(x(\tau))} \sup_{t \in [0,\min\{\theta(z(\cdot)),T_\xi\})} \widetilde{\alpha}(|h(z(t))|)l(t) \\
&\leq \sup_{\widehat{z}(\cdot) \in \mathbf{S}(\xi)} \sup_{t \in [\tau,\min\{\theta(\widehat{z}(\cdot)),\tau+T_\xi\})} \widetilde{\alpha}(|h(\widehat{z}(t))|)l(t-\tau) \\
&\leq \sup_{\widehat{z}(\cdot) \in \mathbf{S}(\xi)} \sup_{t \in [0,\min\{\theta(\widehat{z}(\cdot)),\tau+T_\xi\})} \widetilde{\alpha}(|h(\widehat{z}(t))|)l(t-\tau) \\
&= \sup_{\widehat{z}(\cdot) \in \mathbf{S}(\xi)} \sup_{t \in [0,\min\{\theta(\widehat{z}(\cdot)),\tau+T_\xi\})} \widetilde{\alpha}(|h(\widehat{z}(t))|)l(t)\frac{l(t-\tau)}{l(t)} \\
&\leq \sup_{\widehat{z}(\cdot) \in \mathbf{S}(\xi)} \sup_{t \in [0,\theta(\widehat{z}(\cdot)))} \widetilde{\alpha}(|h(\widehat{z}(t))|)l(t) \cdot \sup_{t \in [0,\tau+T_\xi)} \frac{l(t-\tau)}{l(t)} \\
&\leq V(\xi)\left(\max_{t \in [0,\tau+T_\xi]} \frac{l(t-\tau)}{l(t)}\right). \quad (17)
\end{aligned}$$

We next indicate how $\widetilde{\alpha}$, $\widetilde{\beta}$ and $l$ can be chosen to guarantee that $V$ is an exponential decay RES-Lyapunov function.

Following [32], we make use of the following Lemma on $\mathcal{KL}$ functions ([25], Proposition 7).

**Lemma 6.6** Given a function $\beta \in \mathcal{KL}$ and any number $\lambda > 0$, there exist two functions $\overline{\alpha}_1, \overline{\alpha}_2 \in \mathcal{K}_\infty$ so that
$$\overline{\alpha}_1(\beta(s,t)) \leq \overline{\alpha}_2(s) e^{-\lambda t} \qquad \forall s \geq 0, t \geq 0.$$

With the choice of $\overline{\alpha}_1$ and $\overline{\alpha}_2$ given by Lemma 6.6 for the gain $\beta$ and $\lambda = 2$, the stability condition (13) gives
$$\overline{\alpha}_1(|y(t)|) \leq \overline{\alpha}_2(|\xi|_\omega) e^{-2t} \qquad \forall t \in [0, \theta(x(\cdot)))$$
for each $x(\cdot) \in \mathbf{S}(E)$. Choosing $\widetilde{\alpha} = \overline{\alpha}_1$, we set $l(t) := e^t$, and so can choose the $\mathcal{KL}$ function $\widetilde{\beta}$ as $\widetilde{\beta}(s,t) := \overline{\alpha}_2(s) e^{-t}$.

Finally, we will verify that with this choice of $\widetilde{\alpha}, \widetilde{\beta}$ and $l$, the function $V$ satisfies the exponential decrease requirement (9) with $\alpha_3(r) = r$.

**Proposition 6.7** Given the choice of $\widetilde{\alpha}, \widetilde{\beta}$ and $l$ above, the function $V$ satisfies (9) with $\alpha_3(r) = r$. That is, for each trajectory $x(\cdot) \in \mathbf{S}(E)$,
$$x(s) \in E \text{ for } s \in [0, t] \implies V(x(t)) - V(x(0)) \leq -\int_0^t V(x(s))\, ds.$$

*Proof.* Let $V$ be defined as above, and suppose $x(\cdot) \in \mathbf{S}(E)$ is such that $x(s) \in E$ for $s$ in some interval $[0, t_1]$. Then, as
$$\frac{l(t-\tau)}{l(t)} = e^{-\tau}$$
for all $t \geq 0$, we may take $T_\xi = \infty$ in the decrease argument above (which then holds for all $\tau \in [0, t_1]$), and find that (17) gives
$$V(x(t)) \leq V(x(0)) e^{-t} \qquad \forall t \in [0, t_1]. \tag{18}$$

For $t \in [0, t_1]$, define
$$m(t) := V(x(0)) e^{-t} + \int_0^t V(x(s))\, ds.$$

Then $m(t)$ is an absolutely continuous function whose derivative exists almost everywhere on $[0, t_1]$. We find, from (18),
$$\frac{d}{dt} m(t) = -V(x(0)) e^{-t} + V(x(t)) \leq 0$$
for almost every $t \in [0, t_1]$. Thus $m(t)$ is non-increasing, so in particular, $m(t) \leq m(0)$ for each $t \in [0, t_1]$. That is,
$$V(x(0)) e^{-t} + \int_0^t V(x(s))\, ds \leq V(x(0)) \qquad \forall t \in [0, t_1]. \tag{19}$$

We conclude from (18) and (19) that
$$\begin{aligned} V(x(t)) - V(x(0)) &\leq V(x(0)) e^{-t} - V(x(0)) \\ &\leq -\int_0^t V(x(s))\, ds \end{aligned}$$
for all $t \in [0, t_1]$. ∎

This completes the proof of necessity: we have shown that the function $V$ is lower semicontinuous and satisfies (6) with $\alpha_1(r) = \widetilde{\alpha}(r) l(0)$, $\alpha_2(r) = \widetilde{\beta}(r, 0)$ and (9) with $\alpha_3(r) = r$. ∎

# 7  Discussion

As previously mentioned, the MES property (or more precisely IMES – partial detectability under explicit inputs) is a natural combination of the notions of IOS and IOSS. As such, one would hope that a Lyapunov characterization of the IMES property would include as special cases the existing characterizations for IOS and IOSS (derived in [30] and [14], respectively). The work presented here is a first step toward such a single unifying result.

Several extensions to this result will be needed to complete this program. Firstly, an explicit input can be included by modelling the system as a forced differential inclusion. Secondly, a complete Lyapunov characterization is needed, with no recourse to an additional boundedness assumption. Finally, one would hope to prove that the stability property implies the existence of a smooth Lyapunov function, rather than the discontinuous case described here. When and if these problems are addressed, there will be a single characterization which would encompass the Lyapunov results on ISS, IOS and IOSS.